\documentstyle[twoside]{article}

\def \C{{\Bbb C}}

\def \CP{{\Bbb C}{\Bbb P}}

\def \pii{{\pi_{\infty}}}
\def \R{{\Bbb R}}

\def \O{{\Omega}}

\def \T{{\Theta}}

\def \d{{\partial}}

\def \i{{\sqrt{-1}}}

\def \kk{{\frak k}}
\def \l{{\lambda}}

\def \u{{\frak u}}
\def \su{{\frak s}{\frak u}}
\def \O{{\cal O}}
\def \o{{\omega}}
\def \del{{\partial}}

\def \proof{{\noindent{\it Proof.\ \ }}}
\newtheorem{Th}{THEOREM}[section]

\newtheorem{prop}[Th]{PROPOSITION}

\date{June 25, 2001}

\title{Integrable systems associated with the Bruhat Poisson structures.}
\author{Philip Foth}
\begin{document}
\maketitle
\input amssym.def
\markboth{Philip Foth}{Integrable systems and Bruhat Poisson structures.}

\begin{abstract}
The purpose of this note is to give a simple description of a (complete) 
family of functions in involution on certain hermitian symmetric spaces.
This family obtained via the bi-hamiltonian approach using the Bruhat Poisson
structure is especially simple for projective spaces, where the formulas
in terms of the momentum map coordinates are presented. We show how these 
functions are related to the Gelfand-Tsetlin coordinates. We also show
how the Lenard scheme can be applied.   
\end{abstract}

\section{Introduction.}
\setcounter{equation}{0}
\footnotetext{Research is partially supported by NSF grant DMS-0072520.}

Let $K$ be a compact real form of a complex semi-simple Lie group $G$
and let $H\subset K$ be a subgroup of $K$ defined by $H=K\cap P$, where
$P$ is a parabolic subgroup of $G$ containing a Borel subgroup $B\subset G$. 
The Bruhat Poisson structure $\pii$ on $X=K/H$, first introduced by 
Soibelman \cite{Soi} and Lu-Weinstein \cite{LW}, has the property 
that its symplectic leaves are precisely the Bruhat cells
in $X$. If $T=K\cap B$ is a maximal torus of $K$, then $\pii$ is $T$-invariant. 
Let $\o_s$ (respectively, $\pi_s$) stand for a $K$-invariant symplectic form 
(respectively, dual bi-vector field) on $X$, which we assume now to be a compact 
hermitian symmetric space. It was shown by Khoroshkin-Radul-Rubtsov 
in \cite{KKR} that the two Poisson structures, $\pii$ and $\pi_s$ are 
compatible, meaning that the Schouten bracket of $\pii$ and 
$\pi_s$ vanishes, $[\pii, \pi_s ]=0$. In particular, any bi-vector field of the
form $\alpha \pii +\beta \pi_s$, $(\alpha,\beta )\in \R^2$, is Poisson. 
In this situation, one can
introduce the following family $\{ f_k\}$ of functions: $$f_k:= (\pii^{\wedge k}, 
\o_s^{\wedge k}),$$ obtained by the duality pairing of exterior powers of $\o$ and 
$\pii$. If the (real) dimension of $X$ is equal to $2n$, then we have $n$ functions, 
$f_1$, ..., $f_n$ which may carry some useful information about $X$. 
These functions are in involution with respect to either of the two Poisson 
structures. We make explicit computations for 
$\CP^n$, since this is the only case, where we can present an explicit 
coordinate approach. We show how the function that we have obtained are related 
to the Gelfand-Tsetlin integrable systems studied by Guillemin and Sternberg 
\cite{GSGT}. Analogous statements for other hermitian symmetric spaces will 
appear elsewhere \cite{EF}. In the last part of the paper we make explicit
computations using the Lenard scheme \cite{Magri}. 

\

\noindent{\bf Acknowledgments.} I would like to thank Lu Jiang-Hua and Sam Evens 
for answering many questions regarding Bruhat-Poisson structures. 
Lu Jiang-Hua also provided simple proofs of Propositions 2.1 and 2.2. 
I thank Hermann Flaschka for conversations about integrable systems. I thank 
Yan Soibelman for historical remarks, and Ping Xu for discussions about 
Poisson-Nijenhuis manifolds.  

\

\section{Families of functions in involution.}

Multi-hamiltonian structures are very important in the theory of integrable
systems. Starting with the fundamental works of Magri \cite{Magri}, 
bi- and multi-Hamiltonian structures
found many interesting and fundamental applications, as in \cite{KSM},
\cite{GZ}, \cite{RST} and references therein.   

Let $M$ be a manifold and let $\pi_b$ and $\pi_s$ be two Poisson structures on 
$M$ such that 

\noindent 1. The Poisson structure $\pi_s$ is non-degenerate (so the subscript
$s$ stands for symplectic). 

\noindent 2. The Poisson structures $\pi_s$ and $\pi_b$ are compatible, meaning that 
the Schouten bracket $[\pi_s, \pi_b]$ vanishes. Or, equivalently, for any two 
real numbers $\alpha$ and $\beta$, 
the bi-vector field $\alpha\pi_s + \beta\pi_b$ defines a Poisson structure on $M$. 

\

If ${\rm dim}(M) =2n$, then we can define $n$ functions $f_1$, ..., $f_n$ as 
follows: 

$$ f_j = {\pi_b^j \wedge \pi_s^{n-j}\over \pi_s^n}.$$ 
The operation of division by the top degree bi-vector field makes perfect sense, 
since $\pi_s$ is non-degenerate, and thus in any local coordinate system
$(x_1,$ ... $x_{2n})$ the $2n$-vector field $\pi_s^{n}$ looks like 
$$
\pi_s^{n}=h(x_1, ...,x_{2n}) \del_{x_1}\wedge \cdots\wedge\del_{x_{2n}},
$$
for a non-vanishing function $h(x_1$, ..., $x_{2n})$.  
Equivalently, if $\o_s$ is the symplectic form dual to $\pi_s$, then one can define
$$f_k:= (\pii^{\wedge k}, \o_s^{\wedge k}),$$
where we use the duality pairing 
$$
\Gamma(M, \wedge^{2k}TM)\otimes \Gamma(M, \wedge^{2k}T^*M)\to C^{\infty}(M).
$$

It turns out that this family of functions has the following property.

\begin{prop} The family of functions $f_i$ defined above are in 
involution with respect
to either Poisson structure, $\pi_b$ or $\pi_s$. \end{prop} 

\proof ({\bf J.-H. Lu}) Let $X_i=i_{df_i}\pi_b$ and let $Y_j=i_{df_j}\pi_s$. 
Consider the equality $f_k\pi_s^n = \pi_b^k\wedge\pi_s^{n-k}$ and
compute $L_{X_l}$ of both sides to arrive to the following identity:
$$
{n-k\over k+1} \{ f_{k+1}, f_l\}_s = -\{ f_k, f_l\}_b+nf_k\{ f_1, f_l \}_s,
$$  
where the subscripts $s$ or $b$ indicate with respect to which Poisson 
structure the Poisson bracket is taken. 
Finally, use the induction on $l$. $\bigcirc$

\

\noindent {\it Remark.} 
The approach that we have followed here is intimately related to the
Poisson-Nijenhuis structures, that were studied by Magri and Morosi
\cite{MM}, Kosmann-Schwarzbach and Magri \cite{KSM}, Vaisman \cite{Vais}
and others. The set of our functions $\{ f_j\}$ can be expressed,
polynomially, through the traces of powers of the intertwining operator
corresponding to the Nijenhuis tensor. 

\

Now let us take $M=X$ to be a coadjoint orbit in $\kk^*$, which we assume 
to be a compact hermitian symmetric space. 
We take $\pi_s=\pi$ - the Kirillov-Kostant-Souriau
symplectic structure and $\pi_b=\pii$ - the Bruhat-Poisson structure, which 
is obtained via an identification of $X$ with $K/(P\cap K)$ as in 
Introduction. Under this identification, $\pi$ is $K$-invariant. 
The following was first proved in \cite{KKR}:

\begin{prop} 
If $X$ is a hermitian symmetric space as above, then 
the Poisson structures $\pi$ and $\pii$ are compatible. 
\end{prop}

\proof ({\bf J.-H. Lu}) Let $X$ be a generating vector field for the $K$-action. 
Clearly, the $K$-invariance of $\pi$ implies that $L_X\pi=0$. Since
$\pii$ came from $\kk\wedge\kk$ by applying left and right actions of $K$, 
$L_X\pii$ is obtained from $\delta(X)$ by applying the $K$-action. 
Here, $\delta(X)$ is the co-bracket of $X$, which is an
element of $\kk\otimes \kk$, since we can view $X$ as
an element of $\kk$. Therefore, $L_X\pii$ is a sum of 
wedges of generating vector fields for the action of $K$. Accordingly, 
$$ [L_X\pii, \pi]=0,
$$
which in turn implies that 
$$
L_X[\pii, \pi]=0,
$$
and thus $[\pii, \pi]$ is a $K$-invariant 3-vector field on $X$. 
When $X$ is a hermitian symmetric space, there are none such 
(since the nil-radical of the corresponding parabolic group is 
abelian), so it must be zero. $\bigcirc$

\

Therefore, we have the following

\begin{prop} Let 
$X$ be a coadjoint orbit in $K$. Assume that $X$ is a hermitian symmetric 
space of complex dimension $n$. The above recipe yields $n$ functions 
$(f_1$, ..., $f_n)$ on $X$, which are in involution with respect to either
$\pii$ or $\pi$. 
\end{prop}

The functions $(f_1$, ..., $f_n)$ that we have constructed turn out to be
related to the Gelfand-Tsetlin coordinates in the case when $K=SU(n)$, 
as we will see later on. In the next section we will carry explicit
computations of these functions on the projective spaces.  

\

\section{Computations for the projective spaces.}

Let $\CP^n$ be a complex projective 
space of (complex) dimension $n$, and let $[Z_0:Z_1:...:Z_n]$ be a homogeneous coordinate 
system on it. We use the standard Fubini-Study form $\o$ for $\o_s$ and the following 
description of $\pii$ obtained by Lu Jiang-Hua in \cite{Lu2} and \cite{Lu3}. First, 
we need Lu's coordinates on the largest Bruhat cell, where $Z_0\ne 0$ and we let 
$z_i=Z_i/Z_0$: $$y_i:= {z_i\over\sqrt{1+|z_{i+1}|^2 +\cdots + |z_n|^2}}, \ \ 1\le i\le n.$$

Lu's coordinates are not holomorphic, but convenient for the Bruhat Poisson structure, 
which now assumes the following form
$$\pii = \i \sum_{i=1}^n (1+|y_i|^2)\d_{y_i}\wedge \d_{\bar{y}_i}.$$
In order to be able to compute with $\o$ and $\pii$, we need to move to the polar variables 
$r_i, \phi_j$ defined by $z_i=r_ie^{\i\phi_i}$ and eventually to the momentum map variables 
$x_i, \phi_j$ defined by $$x_i=\delta_{1,i}-{r_i^2\over 1+r_1^2+\cdots r_n^2}.$$ 
These variable are just a slight distortion (for later convenience) of the standard
coordinates on $\R^n$ for the momentum map associated with the maximal compact torus action on 
$\CP^n$. One of the advantages of using this coordinate system is that the Fubini-Study 
symplectic structure has the following simple form: $$\o = \sum_{i=1}^n dx_i\wedge d\phi_i.$$
In fact, the simplest form for the Bruhat Poisson structure is also achieved in this coordinate 
system. 

\begin{prop} The Bruhat Poisson structure $\pii$ on $\CP^n$ can be written in the 
coordinate system $(x_j, \phi_i)$ as $$\pii=\sum_{i=1}^n \T_i\wedge \d_{\phi_i},$$
where $$\Theta_i=(x_1+\cdots +x_i)\d_{x_i} +\sum_{j=i+1}^n x_j\d_{x_j}.$$ 
\end{prop}
\proof The proof of this statement is purely computational. One can introduce
auxiliary variables $q_i=\log(1+|y_i|^2)$, and use those to write 
$$\pii = \sum_{i=1}^n \d_{q_i}\wedge \d_{\phi_i}$$ - the action-angle form for $\pii$. 
Eventually, one can establish the following relations: $x_1=e^{-q_1}$, and for $j>1$,
$$x_j=e^{-(q_1+\cdots +q_j)}-e^{-(q_1+\cdots +q_{j-1})}.$$ The rest is straightforward. 
$\bigcirc$ 

\

Now, one can see that the simple linear and triagonal form of $\pii$ makes 
the computation of the functions $\{ f_i\}$ extremely simple. We will introduce the 
following linear change of variables on $\R^n$:  $$c_k=\sum_{i=1}^k x_i.$$ In these 
variables, the set of functions $\{ f_i\}$ looks as follows.
\begin{Th} The integrals 
$f_i$ (up to constant multiples) arising from the bihamiltonian structure
$(\pi_s, \pii)$ on $\CP^n$ are given by the elementary polynomials in 
$(c_1, ..., c_n)$: 
$$f_1=c_1+\cdots + c_n,$$
$$f_j= \sum_{i_1< ... <i_j}c_{i_1}\cdots c_{i_j},$$
$$f_n=c_1\cdots c_n.$$ \label{Thwe}
\end{Th} 
The explicit nature of these integrals is 
essential in looking at the relation with the certain natural flows \cite{Toda}.
 The hamiltonian $f_1$ in terms of the momentum map variables is given by  
$$f_1= nx_1 +(n-1)x_2 +\cdots + 2x_{n-1}+x_n.$$ Then the
gradient in the momentum simplex has coordinates $\l_i=n+1-i$. Those numbers also
are the 
weights assigned to the vertices (which correspond to the centers of the Bruhat
cells). Thus we arrive to 

\begin{Th} The above flow on $\CP^n$ 
with eigenvalues consecutive integers from $1$ to $n$ 
determines the standard Bruhat 
cell decomposition. \end{Th}

\section{Relation with Gelfand-Tsetlin coordinates.}

When $X=Gr(k)$ - the grassmannian of $k$-planes in $\C^{n+1}$, we have obtained 
$k(n-k+1)$ functions in involution on $X$. Let us recall the standard embedding
$$
\Psi: \ \ F_n \hookrightarrow Gr(1)\times \cdots \times Gr(n),
$$ 
where $F_n$ is the manifold of full flags in $\C^{n+1}$, and the locus of the
embedding is given by the incidence relations. This embedding respects the
KKS Kaehler structures on the manifolds involved, if we would
like to view them as coadjoint orbits in $\kk^*$. Moreover, this embedding is
equivariant with respect to the $K=SU(n+1)$-action. 

Recall the Gelfand-Tsetlin system on $F_n$. We fix the orbit type of $F_n$,
i.e. we fix the eigenvalues $\i \l_i$ and order them, so 
$\l_1 > \l_2 > \cdots > \l_{n+1}$. For convinience and easier visualization, 
we will assume that $\l_{n+1}=0$ (so all the eigenvalues are
non-negative), which will correspond to working with $\u_{n+1}^*$ rather 
than with $\su_{n+1}^*$. For convenience, we also identify the Lie algebra
$\u_k$ with its dual via $-{\rm Tr}(AB)$. The Gelfand-Tsetlin system looks like
\cite{GS5}, \cite{GSGT}:

$$
\begin{array}{ccccccccc}
\l_1 & > & \l_2 & > & \cdots & > & \l_n & > & 0 \\
{} & \mu_1^1 & \ge & \mu_2^1 & \ge & \cdots & \ge & \mu_n^1 & {}\\

{} & {} & {} & \cdots & \cdots & \cdots & {} & {} & {} \\

{} & {} & {} & \mu_1^{n-1} & \ge & \mu_2^{n-1} & {} & {} & {} \\
{} & {} & {} & {} & \mu_1^n & {} & {} & {} & {} 
\end{array},
$$
where the $i$-th row corresponds to the projection
$\u_{n+1}^*\to \u_{n+2-i}^*$, which is dual to the embedding
$U(n+2-i)\hookrightarrow U(n+1)$ in the left upper corner. 
The eigenvalues $\mu_i^j$ together with $\l_i$'s satisfy
the interlacing property. 

The picture above can be adapted to any orbit, in
particular to $Gr(k)$, where we would take
$\l_1=\cdots =\l_k>0$, and other $\l$'s equal to zero. 
When $k$ varies from $1$ to $n$, the picture above acquires more
and more non-zero elements. At each step, while going from
level $k$ to $k+1$, we will get new integrals on $Gr(k+1)$, 
which we can pull back to $F_n$ using $\Psi$. 

Our goal is to relate the integrals $f_j$, that we obtained in 
Section 3 using the bi-hamiltonian approach on hermitian symmetric 
spaces, and the Gelfand-Tsetlin coordinates. We will start working 
with $M = \CP^n$, the complex projective space.   

Let $B$ be the $(n+1)\times (n+1)$ matrix, representing an element of 
$\u(n+1)^*$ such that the only non-zero element of $B$ is $\i\l$, 
located in the very left upper place. The coadjoint orbit $\O_B$ of 
$B$ is isomorphic to $\CP^n$, where the identification goes as follows. 
Any element in the coadjoint orbit of $B$ can be viewed as $ABA^{-1}$, 
where $A\in U(n+1)$. Let $(a_{ij})$ be the entries of $A$. Then 
the identification 
$$
w: \ \O_B\to \CP^n
$$ is given by 
$$
w(ABA^{-1}) = [a_{11}: a_{21}: \cdots : a_{n+1,1}],
$$ in terms of a homogeneous coordinate system 
$[Z_0: \cdots : Z_n]$ on $\CP^n$.
We suspect that the following is well-known, and in any case, is not
hard to compute, that the Gelfand-Tsetlin coordinates are:
$$
\mu_r^k=0 \ \ {\rm for} \ \ r\ne 1, 
$$ $$ 
\mu_1^{k} = \lambda (x_1+ \cdots + x_{n-k+1}), 
$$           
where $(x_1, ..., x_n)$ are the momentum map coordinates that we used in the
previous section. We arrive to the conclusion that the Gelfand-Tsetlin
coordinates $\{ \mu_1^k\}$ coincide 
(up to the multiple of $\lambda$, which we can assume equal to one)
with the coordinates $\{ c_k \}$ introduced in the previous section.  
Now, it remains to notice that the Theorem \ref{Thwe}
from the previous Section
immediately yields

\begin{Th} The complete family of integrals in involution $\{ f_i \}$ 
on $\CP^n$ obtained using the bi-hamiltonian
approach with respect to the Bruhat Poisson structure and an invariant
symplectic structure are expressed by the elementary polynomials in the 
Gelfand-Tsetlin coordinates. 
\end{Th}

We prove a similar result for other hermitian symmetric spaces in a 
forthcoming paper \cite{EF}.

\section{Comparison to the Lenard scheme.}

Recall the following result \cite{Magri}. If $\alpha\pi_0 + \beta\pi_1$
is a Poisson pencil on a manifold $M$, and $V$ a vector field, preserving
this pencil, then there exists a sequence of smooth functions $\{ g_i\}$
on $M$, such that $g_1$ is the Hamiltonian of $V$ with respect to $\pi_0$
and the vector field of the $\pi_0$-hamiltonian $f_j$ is the same as the
vector field of the $\pi_1$-hamiltonian $f_{j+1}$: 
$$i_{df_j}\pi_0 = i_{df_{j+1}}\pi_1.$$
Moreover, the functions in the family $\{ f_j\}$ are in involution with 
respect to both $\pi_0$ and $\pi_1$. 

Our goal in this section is to show that if we start with $M=\CP^n$, and 
take the pencil $(\pi_s, \pii)$ as before, then 
there is a natural choice of $V$ on $\CP^n$ leading to a completely 
integrable systems, and the integrals $\{ g_j \}$ in question can be easily 
expressed in terms of the coordinates $(c_1, ..., c_n)$ that we introduced in 
Section 3.  

It is a matter of a simple computation that if we start with a hamiltonian 
$g_1 = a_1 x_1+\cdots +a_n x_n$, where $(x_1, ..., x_n)$ are the momentum map
coordinates as before, then the corresponding initial vector field
$V$ is given by 
$$
V = i_{dg_1}\pii = \sum_j [a_j(x_1+\cdots x_j)+a_{j+1}x_{j+1}+\cdots + 
a_nx_n]\partial_{\phi_j}.
$$
From this, one can compute $\displaystyle{g_2=\sum_j {a_j\over 2}x_j^2 +
\sum_{l < k} a_kx_lx_k},$ etc. 
An interesting choice for $g_1$ turns out to be 
$$
g_1=c_1+\cdots + c_n = nx_1+(n-1)x_2 +\cdots + x_n,
$$
which coincides with $f_1$ from Section 3. The reason for this choice is 

\begin{prop} The Lenard scheme associated to the Poisson pencil 
$(\pii$, $\pi_s)$ on $\CP^n$ which starts with $g_1=c_1+\cdots + c_n$ and 
$$
V=\sum_j[(n-j+1)(x_1+\cdots + x_j)+(n-j)x_{j+1}+\cdots + 2x_{n-1} +
x_n]\partial_{\phi_j},
$$
yields
$$
g_k=c_1^k+c_2^k+\cdots + c_n^k,
$$
which determines a completely integrable bi-hamiltonian system on $\CP^n$. 
\end{prop}
\proof With all the explicit formulas that we have presented in this paper, 
the proof is a simple computation. $\bigcirc$

\

We should remark, that the constants $(a_1, ...,a_n)$ for the first hamiltonian
in the Lenard scheme have to be chosen with care for two reasons. First, the 
computations are not simple for an arbitrary choice. Second, as the next
example shows, we do not always arrive to a completely integrable system. 

\

\noindent{\bf Example.} If one takes $g_1=x_1+\cdots + x_n$, and 
$V=\sum_j (x_1+\cdots + x_n)\partial_{\phi_j}$, then applying the above 
scheme, one would obtain 
$$  
g_k=(x_1+\cdots + x_n)^k = (g_1)^k. 
$$ 
The differentials of all functions in this family are clearly
linearly dependent.

\vskip 0.3in

Department of Mathematics 
\\ University of Arizona
\\ Tucson, AZ 85721-0089
\\ foth@math.arizona.edu

\

\noindent{\it AMS subj. class.}: \ \ primary 58F07, secondary 53B35, 53C35.


\thebibliography{123}

\bibitem{EF}{S. Evens and P. Foth. {\it In preparation.}} 

\bibitem{GZ}{I. M. Gelfand and I. Zakharevich. Webs, Lenard schemes, and
the local geometry of bi-Hamiltonian Toda and Lax structures. {\it Selecta
Math. (N.S.)}, {\bf 6}: 131-183, 2000.}

\bibitem{GS5}{V. Guillemin and S. Sternberg. On collective complete
integrability according to the method of Thimm. {\it Ergod. Th. 
Dynam. Sys.}, {\bf 3}: 219-230, 1983.}

\bibitem{GSGT}{V. Guillemin and S. Sternberg. The Gelfand-Tsetlin system
and quantization of the Complex Flag Manifolds. {\it J. Funct. An.}, 
{\bf 52}: 106-128, 1983.} 

\bibitem{KKR}{S. Khoroshkin, A. Radul, and  V. Rubtsov. A family of Poisson 
structures on  Hermitian symmetric spaces. 
{\it Comm. Math. Phys.}, {\bf 152}: 299-315, 1993.}

\bibitem{Toda}{J. Moser. Various aspects of integrable Hamiltonian systems. 
in {\it Dynamical systems, Progress Math.}, {\bf 8}: pp.233-289, Birkhauser, Boston, 
1980.}

\bibitem{KSM}{Y. Kosmann-Schwarzbach and F. Magri. 
Poisson-Nijenhuis structures.
{\it Ann. Inst. H. Poincare Phys. Theor.}, {\bf 53}: 35-81, 1990.} 

\bibitem{LW}{J.-H. Lu and A. Weinstein. Poisson-Lie groups, dressing transformations, 
and Bruhat decompositions. {\it J. Diff. Geom.}, {\bf 31}: 501-526, 1990.}

\bibitem{Lu2}{J.-H. Lu. Coordinates on Shubert cells, Kostant's harmonic forms,and 
the Bruhat Poisson structures on $G/B$. {\it Transform. groups}, 
{\bf 4}: 355-374, 1999.}

\bibitem{Lu3}{J.-H. Lu. Darboux coordinates for the Bruhat Poisson structure
on $\CP^n$. {\it Preprint.}}

\bibitem{Magri}{F. Magri. A simple model of the integrable Hamiltonian equation.
{\it J. Math. Phys.}, {\bf 19}: 1156-1162, 1978.}

\bibitem{MM}{F. Magri and C. Morosi. A geometrical characterization of
integrable Hamiltonian systems through the theory of Poisson-Nijenhuis
manifolds. {Quaderno S.}, {\bf 19}, Univ. Milan, 1984.} 

\bibitem{RST}{A. G. Reyman and M. A. Semenov-Tian-Shansky. Group-theoretical
methods in the theory of finite-dimensional integrable systems. In {\it
Dynamical systems. VII. Integrable systems, nonholonomic dynamical systems.}
Encyclopaedia of Mathematical sciences. {\bf 16}: 116-226, Springer-Verlag, 
Berlin, 1994.}

\bibitem{Soi}{Y. Soibelman. The algebra of functions on a compact quantum group 
and its representations. {\it Leningrad J. Math.}, {\bf 2}: 161-178, 1991.}

\bibitem{Vais}{I. Vaisman. Complementary 2-forms of POisson structures. 
{\it Composito Math.}, {\bf 101}: 55-75, 1996.}
\end{document}